\documentclass{sigma_arx}

\renewcommand{\phi}{\varphi}

\newcommand{\C}{\mathbb C}
\newcommand{\R}{\mathbb R}
\newcommand{\Z}{\mathbb Z}

\newcommand{\vect}{\mathcal V}
\newcommand{\diff}{\mathrm{Diff}}
\newcommand{\Aut}{\mathrm{Aut}}
\newcommand{\aut}{\mathfrak{aut}}

\newcommand{\gau}{\mathfrak{gau}}

\newcommand{\Hom}{\mathrm{Hom}}
\newcommand{\End}{\mathrm{End}}

\newcommand{\Hol}{\mathrm{Hol}}
\newcommand{\flhoBun}{\mathrm{Bun}^{\flat}_{\mathrm{hol}}}

\newcommand{\ZDR}{Z^1_{\mathrm{dR}}}
\newcommand{\ZDRhol}{Z^1_{\mathrm{dR,h}}}
\newcommand{\HDR}{H^1_{\mathrm{dR}}}
\newcommand{\HDRhol}{H^1_{\mathrm{dR,h}}}

\newcommand{\g}{\mathfrak g}

\newcommand{\GS}{G ^{\Sigma}}
\newcommand{\Ad}{\mathrm{Ad}}

\newcommand{\id}{\mathrm{id}}

\begin{document}

\renewcommand{\PaperNumber}{001}

\FirstPageHeading

\ShortArticleName{Complex Structures on Principal Bundles}

\ArticleName{Complex Structures on Principal Bundles}

\Author{Martin LAUBINGER}

\AuthorNameForHeading{M. Laubinger}

\Address{Louisiana State University, Baton Rouge, Louisiana, USA} 
\EmailD{martinl@math.lsu.edu} 
\URLaddressD{http://www.math.lsu.edu/\~{}martinl/} 


\Abstract{
Holomorphic principal $G$- bundles over a complex manifold $M$ can be studied using non-abelian cohomology groups $H^1(M,G).$ On the other hand, if $M=\Sigma$ is a closed Riemann surface, there is a correspondence between holomorphic principal $G$-bundles over $\Sigma$ and coadjoint orbits in the dual of a central extension of the Lie algebra $C^\infty(\Sigma, \g).$ We review these results and provide the details of an integrability condition for almost complex structures on smoothly trivial bundles. This article is a shortened version of the author's Diplom thesis.}

\Keywords{holomorphic principal bundles; almost complex structures; non-abelian cohomology; current group; current algebra; coadjoint orbits; central extensions; Riemann surface}

\Classification{58D27; 81R10; 32Q60}

%
%
%
%
%
%
\section{Introduction}
For a Lie algebra valued 1-form $\omega,$ one can define a 2-form $F(\omega)=d\omega+ [\omega, \omega],$ which plays an important role in differential geometry. For instance if $\omega$ is a connection 1-form in a principal $G$-bundle, then $F(\omega)=0$ means that the corresponding connection is flat. 
On complex manifolds, the de Rham differential decomposes as $d = \partial + \bar \partial,$ and we define $\bar F(\omega) = \bar \partial \omega = [ \omega, \omega].$ This expression is used in \cite{Onishchik} in a theorem which states that the equation $f^{-1}.\bar \partial f = \omega$ is locally solvable for $f$ if and only if $\bar F(\omega)=0.$ One of the main results of the present paper, Theorem \ref{MainInt}, is concerned with certain almost complex structures $I_\omega$ on products $M \times G$ associated with a $(0,1)$-form $\omega.$  We show that the structure $I_\omega$ is integrable, thus yielding a holomorphic principal bundle, if and only if $\bar F(\omega)=0.$

The present article is based on the author's Diplom thesis \cite{Laubinger}. The starting point of our thesis was section 3 in \cite{EtingofFrenkel}, where holomorphic principal bundles over a Riemann surface $\Sigma$ are related to coadjoint orbits of extensions of the so-called current algebra $C^\infty( \Sigma, \g).$ We give a proof of this correspondence in Theorem \ref{MainThm}, which provides more details than the discussion in \cite{EtingofFrenkel}. Also, we review some of the results from \cite{Onishchik} on classification of holomorphic principal bundles via non-abelian cohomology.

Let us summarize the content of the following four sections. Section 2 contains a summary of the basic facts about almost complex and complex manifolds that are needed later. We formulate an integrability condition for almost complex structures, and define the Maurer-Cartan equation $F(\omega)=0$ as well as the similar equation $\bar F(\omega)=0,$ which we call the Holomorphic Maurer-Cartan equation.

In Section 3, we review material about principal $G$-bundles. We define and characterize flat bundles, and define the gauge group of a principal bundle. In the case of a trivial bundle, the gauge group can be identified with the group $C^\infty(M,G)$ of smooth $G$-valued maps. The latter group is also frequently called \emph{current group}, and its Lie algebra $C^\infty(M, \g)$ is the \emph{current algebra}. In Section 3.5 we focus on bundles over a closed Riemann surface. In case of a surface of genus 1 and a compact group $G=K,$ we use a result of Armand Borel on automorphisms of compact Lie groups to describe the space of equivalence classes of principal $G$-bundles quite explicitly. 

Sections 4 and 5 constitute the main part of the present article. In Section 4 we describe almost complex structures on the total space $P$ of a smoothly trivial principal $G$-bundle in terms of differential forms of type $(0,1)$ on the base manifold. We give criteria for the integrability of the almost complex structure, as well as for the holomorphic equivalence of two such structures. 

In Section 5 we define for each holomorphic 1-form $\eta$ a central extension $E_\eta$ of $C^\infty(\Sigma, \g).$ Here we assume that $\g$ is simple, so that there is an invariant inner product $\langle \cdot, \cdot \rangle.$ The extension is defined by the cocycle $\Omega(f,g) = \int_\Sigma \eta \wedge \langle f, dg \rangle.$ Then we define a pairing between the extension $E_\eta$ and the vector space $E^*=\C \oplus \Omega^{(0,1)}(\Sigma, \g).$ This pairing establishes the connection to holomorphic principal bundles, since each element of $\Omega^{(0,1)}(\Sigma, \g)$ yields a holomorphic bundle structure on the product $\Sigma \times G.$ Finally we can state the main result Theorem \ref{MainThm}, which relates coadjoint orbits in $E^*$ to equivalence classes of holomorphic principal bundles.

%
\section{Almost Complex Structures}

We assume familiarity with the basic constructions of differential geometry of manifolds and vector bundles. The main references for complex and almost complex manifolds are \cite{KobayashiNomizu:1963}, Chapter IX, and \cite{Wells}.

\begin{definition}
Let $M$ be a real smooth $m$-dimensional manifold. An \emph{almost complex structure on $M$} is a smooth section $I$ of $\End(TM)$ such that $I^2_x = - \id_{T_xM}$ for each $x \in M.$ The pair $(M,I)$ is an \emph{almost complex manifold}. Note that $I^2_x = - \id_{T_xM}$ is only possible if the dimension of $M$ is an even number $m = 2n.$ If $(M,I)$ and $(N,J)$ are almost complex manifolds and $f:M \rightarrow N$ a smooth map, then we say that $f$ is \emph{holomorphic} if its differential intertwines the complex structures, i.e. if $df \circ I = J \circ df.$
\end{definition}

\begin{example} \label{ex1}
Every complex $n$-dimensional manifold also carries the structure of a $2n$-dimensional real smooth manifold. Its tangent spaces are complex vector spaces, and if we define $I_x$ to be multiplication with $i$ in $T_xX,$ then $(M,I)$ is a $2n$-dimensional almost complex manifold.
\end{example}

\begin{definition}
We call an almost complex structure $I$ on a real smooth manifold $M$ \emph{integrable} if $M$ admits the structure of a complex manifold such that $I_x$ corresponds to multiplication by $i$ as in Example \ref{ex1}.
\end{definition}

Not every even dimensional real smooth manifold admits an almost complex structure and not every almost complex manifold is integrable. The four-sphere $S^4$ for example does not allow an almost complex structure. One of the main theorems of this article is concerned with the integrability of an almost complex structure on the total space of a principal bundle.

Note that if $I$ is an almost complex structure on $M,$ then an eigenvalue $\lambda$ of $I_x$ has to satisfy $\lambda^2 = -1,$ hence $I_x$ has no real eigenvalues. However, if we extend $I_x$ to the complexification $(T_xM)_\C$ then it has eigenspaces of complex dimension $n$ corresponding to the eigenvalues $i$ and $-i.$ This leads to the following definition.

\begin{definition} \label{def:Decomp}
Let $(M,I)$ be an almost complex manifold. Then $I_x$ extends to an endomorphism of $T_xM_\C$ and also of the dual $T_xM_\C^*$ via $I_x\xi(v)=\xi(I_xv)$ for $\xi \in T_xM_\C^*$ and $v \in T_xM_\C.$ Let 
\[
TM_\C = TM^{(1,0)}\oplus TM^{(0,1)}, \quad TM_\C^* = TM^*_{(1,0)} \oplus TM^*_{(0,1)}.
\]
be the corresponding decompositions in subbundles consisting of $i$-eigenvectors and $(-i)$-eigenvectors, respectively. We then call sections of $TM^{(1,0)}$ and $TM^{(0,1)}$ \emph{vector fields of type $(1,0)$ and $(0,1)$} respectively, and write $\vect(M)=\vect^{(1,0)}(M) \oplus \vect^{(0,1)}(M).$ The exterior products of $TM_\C^*$ decompose into a direct sum with summands of the form
\[
\sideset{}{^{(p,q)}} \bigwedge TM_\C^* = \text{span}_\C \left\{ \alpha \wedge \beta \, | \, \alpha \in \wedge^p TM_{(1,0)}^*, \beta \in \wedge^q TM_{(0,1)}^* \right\}.
\]
Smooth sections of these summands are called \emph{differential forms of type (p,q),} and $\Omega^{(p,q)}(M)$ is the vector space of forms of type $(p,q).´$ Usually one chooses local coordinates $x_1, \dots, x_n, y_1, \dots, y_n$ such that the vectors
$\frac{\partial}{\partial z_k} := \frac{1}{2}\left( \frac{\partial}{\partial  x_k} -i\frac{\partial}{\partial y_k}\right)$ and 
$\frac{\partial}{\partial \bar z_k} := \frac{1}{2} \left(\frac{\partial}{\partial x_k} +i\frac{\partial}{\partial y_k} \right)$ 
form local frames for $T^{(1,0)}M$ and for $T^{(0,1)}M,$ respectively. Then a form $\omega$ of type $(p,q)$ can be expressed locally as
\[
\omega|_U = \sum_{i_0 < \dots i_{p+q}} \alpha_{i_0, \dots, i_{p+q}} dz_{i_0} \wedge \dots \wedge dz_{i_p} \wedge d \bar z_{i_{p+1}} \wedge d \bar z_{i_{p+q}}
\]
with smooth coefficients $\alpha_{i_0, \dots, i_{p+q}}.$
\end{definition}

\begin{definition}
If $(M,I)$ is an almost complex manifold, then the torsion $N_I(X,Y)$ of two vector fields $X,Y \in \vect(M)$ is defined to be the vector field given by
\[
N_I(X,Y)= \left[ IX,IY \right] -\left[ X,Y\right] -I \left[IX,Y \right] -I\left[X,IY \right].
\]
\end{definition}
The following theorem characterizes integrable almost complex manifolds.

\begin{theorem} \label{integrability}
For an almost complex manifold $(M,I),$ the following are equivalent:
\begin{itemize}
\itemsep=0pt
\item $I$ is integrable
\item The decomposition $\vect(M)=\vect^{(1,0)}(M) \oplus \vect^{(0,1)}(M)$ is stable under formation of commutators. In other words, $\vect^{(1,0)}(M)$ and $\vect^{(0,1)}(M)$ are Lie subalgebras of $\vect(M).$
\item $d \Omega^{(p,q)}(M) \subset \Omega^{(p+1,q)}(M) \oplus \Omega^{(p,q+1)}(M)$
\item The torsion $N_I$ vanishes identically
\end{itemize}
\end{theorem}
\begin{proof}
See \cite{KobayashiNomizu:1963}, Chapter IX, Theorem 2.8.
\end{proof}

\begin{definition}
By Theorem \ref{integrability}, if $M$ is a complex manifold, then the de Rham differential can be written as $d = \partial + \bar \partial,$ where $\partial \omega \in \Omega^{(p+1,q)}(M)$ and $\bar \partial \omega \in \Omega^{(p,q+1)}(M)$ for a form $\omega$ of type $(p,q).$ We define the vector space of \emph{holomorphic $k$-forms} to be 
\[
\Omega^k_h(M) = \{ \omega \in \Omega^{(k,0)}(M) \; | \; \bar \partial \omega = 0 \}.
\]
The holomorphic 0-forms are exactly the holomorphic functions on $M.$
\end{definition}

\begin{theorem}
Let $M$ be a complex $m$-dimensional manifold. Then we have an exact sequence of vector spaces
\[
0 \rightarrow  \Omega^p_h(M)  \stackrel{i}{\rightarrow}  \Omega^{(p,0)}(M) \stackrel{\bar \partial}{\rightarrow}  \Omega^{(p,1)}(M) \stackrel{\bar \partial}{\rightarrow}  \cdots  \stackrel{\bar \partial}{\rightarrow}  \Omega^{(p,m-p)}(M) \rightarrow 0 
\]
for each $p <m$ and a complex of vector spaces
\[
0 \rightarrow \C   \stackrel{i}{\rightarrow}  \Hol(M) \stackrel{d}{\rightarrow} \Omega^1_h(M)  \stackrel{d}{\rightarrow} \Omega^2_h(M) \stackrel{d}{\rightarrow} \cdots \stackrel{d}{\rightarrow} \Omega^m_h(M)  \rightarrow 0.
\]
\end{theorem}
\begin{proof}
See \cite{Wells}, Chapter II, Examples 2.12 and 2.13.
\end{proof}

\begin{definition}
If $(M,I)$ is an almost complex manifold and $V$ a complex vector space, then we let $\Omega^p(M,V)=\Omega^p(M) \otimes_\C V$ be the vector space of \emph{$p$-forms on $M$ with values in $V.$} This can be done in particular when $V=\g$ is a complex Lie algebra, and in this case we also define
\[
\left[ \cdot, \cdot \right]: \Omega^1(M,\g) \times \Omega^1(M,\g) \rightarrow \Omega^2(M,\g)
\]
by
\[
\left[ \alpha, \beta \right]_p(X,Y) := \left[ \alpha_p(X),\beta_p(Y) \right]-\left[ \alpha_p(Y), \beta_p(X) \right].
\] 
\end{definition}

\begin{remark} \label{rem:typeTwozero}
It is not hard to show that if $\omega \in \Omega^{(0,1)}(M,\g),$ then $\left[ \omega, \omega \right] \in \Omega^{(0,2)}(M,\g).$ Also note that $[\omega, \eta] = [\eta, \omega]$ for 1-forms $\omega$ and $\eta.$
\end{remark}

\begin{definition} \label{def:MCeqn}
If $\omega$ is a Lie algebra valued 1-form on a manifold $M,$ then we let
\[
F(\omega) = d\omega + \frac{1}{2} \left[ \omega, \omega \right],
\]
and if $M$ is an almost complex manifold, we let
\[
\bar F(\omega) = \bar \partial \omega + \frac{1}{2} \left[ \omega, \omega \right].
\]
Note that by Lemma \ref{rem:typeTwozero}, if $\omega \in \Omega^{(0,1)}(M, \g),$ then $\bar F(\omega) \in  \Omega^{(0,2)}(M,\g).$ The equation $F(\omega)=0$ is called the \emph{Maurer-Cartan equation, or MC equation.} We call the equation $\bar F(\omega)=0$ \emph{holomorphic Maurer-Cartan equation, or HMC equation.}
\end{definition}
%
\section{Principal $G$- Bundles}

We recall the basic definitions concerning principal $G$-bundles (see for instance \cite{Sharpe}).
\begin{definition}
Let $M$ be a smooth manifold and $G$ a Lie group. A \emph{principal $G$-bundle} is a smooth principal fiber bundle $q:P \rightarrow M$ with typical fiber $G$ and a smooth right action of $G$ on $P,$ such that the bundle charts $\psi_U: q^{-1}(U) \rightarrow U \times G$ are equivariant with respect to the right $G$-action $(x,h)g =(x,hg)$ on $U \times G.$
\end{definition}

\begin{definition}
Principal $G$-bundles can be described by \v Cech-2-cocycles as follows. If $\{ U_j \}$ is an open cover of $M,$ we use the convention $U_{jk}= U_j \cap U_k.$ Furthermore, let $P_U= q^{-1}(U).$ If $\psi_j:P_{U_j} \rightarrow U_j \times G$ are the bundle charts of a principal $G$-bundle $P,$ then the change of coordinates $\psi_j \circ \psi_k^{-1}:U_{jk} \times G \rightarrow U_{jk} \times G$ can be written as $(x,h) \mapsto (x, g_{jk}(x)h),$ where the $g_{jk}:U_{jk} \rightarrow G$ are smooth maps. We call the bundle \emph{holomorphic} if $M$ is a complex manifold, $G$ is a complex Lie group and the change of coordinates $\psi_j \circ \psi_k^{-1}$ or equivalently the $g_{jk}$ are holomorphic maps.
\end{definition}

\begin{lemma}
For a principal $G$-bundle $q:P \rightarrow M$ let us denote the right action of $G$ by  $\sigma:  P \times G \rightarrow P,$ and let $\vect(P)$ denote the Lie algebra of vector fields on $P.$ Then the map
\[
\dot{\sigma}: \g \rightarrow \vect(P), \quad \dot{\sigma}(x)_p := d\sigma_{(p,e)}(0,x)
\]
is a homomorphism of Lie algebras.
\end{lemma}

\begin{definition}
Let $\rho_g: P \rightarrow P, p \mapsto p.g = \sigma(p,g)$ be the right action of $g$ on $P.$ Then $\omega \in \Omega^1(P,\g)$ is called \emph{connection 1-form} if it satisfies the following conditions:
\begin{enumerate}
\item $\omega_p(\dot{\sigma}(x)_p)=x$ for all $x \in \g$ and $p \in P.$
\item $\omega \circ d\rho_g = \Ad(g^{-1}) \circ \omega.$
\end{enumerate}
The \emph{curvature} of a connection 1-form is the 2-form $F(\omega) = d\omega + \frac{1}{2} [ \omega, \omega ] \in \Omega^2(P,\g).$ If the curvature vanishes, that is, if $\omega$ satisfies the MC equation, then the connection is called \emph{flat.} A principal $G$-bundle is called \emph{flat} if there exists a flat connection 1-form on $P.$
\end{definition}

\begin{definition}
Let $\rho: \pi_1(M) \rightarrow G$ be a homomorphism and $\widetilde M$ the universal covering space of $M.$ The fundamental group $\pi_1(M)$ acts on $\widetilde M$ from the left by covering transformations, and then on $\widetilde M \times G$ via $\gamma.(m,g)=(\gamma.m,\rho(\gamma)g).$ It can be shown that the quotient of $\widetilde M \times G$ by this action is a manifold, and the induced projection onto $M$ defines a principal $G$-bundle which we call $P_\rho.$
\end{definition}

\begin{lemma}
For a principal $G$-bundle $P,$ the following are equivalent:
\begin{itemize}
\item $P$ admits a flat connection
\item $M$ admits a cover $\{U_i\}$ such that $P$ can be defined by a \v Cech cocycle consisting of constant maps $g_{ij}.$
\item $P$ is isomorphic to $P_\rho$ for a homomorphism $\rho: \pi_1(M) \rightarrow G.$
\end{itemize}
\end{lemma}

\begin{definition}
If $P$ is a principal $G$-bundle, let $\Aut(P)$ denote the group of automorphisms of $P,$ that is diffeomorphisms $\phi$ of $P$ which satisfy $\phi(pg)=\phi(p)g$ for all $g \in G.$ There is a natural homomorphism of $\Aut(P)$ onto $\diff(M),$ and the kernel of this homomorphism is the \emph{gauge group}. 
Similarly, $\aut(P)$ is defined to be the Lie algebra of invariant vector fields on $P,$ that is, vector fields $\xi$ which satisfy $\rho_g^* \xi = \xi$ for all $g \in G.$ The kernel of the projection $\aut(P) \rightarrow \vect(M)$ is thought of as the Lie algebra of the gauge group. Connections in $P$ can be identified with linear splittings of 
\[
0 \rightarrow \gau(P) \rightarrow \aut(P) \rightarrow \vect(M) \rightarrow 0
\]
and the curvature of  a connection measures the failure of such a splitting to be a homomorphism of Lie algebras.
\end{definition}
\subsection{Fundamental Theorem}

The main purpose of this subsection is to define the period homomorphism $\mathrm{per}_\omega: \pi_1(M) \rightarrow G$ associated to a solution $\omega$ of the MC-equation. See \cite{Sharpe} for proofs of Theorems \ref{fundThmSC} and \ref{fundThm}.

\begin{definition}
If $G$ is a Lie group, let $\lambda_g$ and $\rho_g$ denote left and right translation by $g \in G.$ We get corresponding actions of $G$ on $TG$ as follows. If $\xi \in T_hG,$ let $g. \xi= d_h\lambda_g(\xi) \in T_{gh}G$ and $\xi. g= d_h\rho_g(\xi) \in T_{hg}G.$ Now if $f \in C^\infty(M,G),$ we define the \emph{left logarithmic derivative of $f$} to be the $\g$-valued 1-form on $M$ given by $\delta(f)_m= f(m)^{-1}.df_m: T_mM \rightarrow \g.$ For short, we write $\delta(f) = f^{-1}.df,$ and we define the \emph{right logarithmic derivative} by $\delta^r(f) = df.f^{-1}.$ 
\end{definition}

\begin{theorem}[Fundamental Theorem for simply connected manifolds] \label{fundThmSC}
Let $M$ be a simply connected manifold, $G$ a Lie group and $\omega \in \Omega^1(M,\g).$ Then there exists an 
$f \in C^{\infty}(M,G)$ with $\delta(f)=\omega$ if and only if $\omega$ satisfies the MC equation. For two elements $f_1,f_2 \in C^{\infty}(M,G)$ with $\delta(f_1)= \delta(f_2)$ there exists an element $g \in G$ with $f_2 = \lambda_g \circ f_1.$
\end{theorem}

\begin{lemma}
If $q: \widetilde M \rightarrow M$ is the universal  cover of a smooth manifold $M,$ let $q^* \omega$ denote the pull-back of a $\g$-valued 1-form $\omega$ on $M$ to $\widetilde M.$ Let $\mu_\gamma$ denote the left action of an element $\gamma \in \pi_1(M)$ on $\widetilde M.$ Then if $f \in C^\infty(\widetilde M, G)$ satisfies $\delta(f) = q^* \omega,$ so does $f \circ \mu_\gamma.$ In particular, by the Fundamental Theorem, $f \circ \mu_\gamma = \lambda_{c(\gamma)} \circ f$ for some $c(\gamma) \in G.$ The map $\gamma \mapsto c(\gamma)$ defines a group homomorphsim $c: \pi_1(M) \rightarrow G.$
\end{lemma}

\begin{definition}
Let the situation be as in above lemma, and let $\omega \in \Omega^1(M,\g)$ satisfy the MC equation. Fix $m_0 \in \widetilde M$ and let $S(\omega) \in C^\infty( \widetilde M, G)$ be the unique function satisfying $\delta(S(\omega)) = q^* \omega$ and $S(\omega)(m_0)=e.$ Then the corresponding homomorphism $c:\pi_1(M) \rightarrow G$ is called the \emph{period homomorphism of $\omega$} and is denoted $\text{per}_\omega.$
\end{definition}

\begin{theorem}[Fundamental Theorem] \label{fundThm}
Let $M$ be a connected manifold, $G$ a Lie group and $\omega \in \Omega^1(M,\g).$ Then there is an 
$f \in C^\infty(M,G)$ such that $\omega = \delta(f)$ if and only if $\omega$ satisfies the MC equation
and its period homomorphism $\mathrm{per}_\omega$ is constant.
For two elements $f_1,f_2 \in C^\infty(M,G)$ with $\delta(f_1)=\delta(f_2)$ there exists an element 
$g \in G$ with $f_2 = f_1 \circ \lambda_g.$
\end{theorem}

The following holomorphic version of the fundamental theorem holds:
\begin{theorem}
Let $\omega \in \Omega_h^1(M,\g).$ Then there is an 
$f \in \Hol(M,G)$ such that $\omega = \delta(f)$ if and only if $\omega$ satisfies the MC equation
and its period homomorphism $\mathrm{per}_\omega$ is constant.
For two elements $f_1,f_2 \in \Hol(M,G)$ with $\delta(f_1)=\delta(f_2)$ there exists an element $g \in G$ with
$f_2 = f_1 \circ \lambda_g.$
\end{theorem}
\begin{proof}
Given $\omega \in \Omega^1_h(M,\g),$ the Fundamental Theroem \ref{fundThm} assures us that there is a map
$f \in C^\infty(M,G)$ with $\delta(f)=\omega$ if and only if $\omega$ satisfies the MC equation and 
the period homomorphism $\mathrm{per}_\omega$ is trivial.
It remains to show that $\delta(f)=\omega$ implies that $f$ is holomorphic. 
But as $\delta(f)= f^{-1}.df$ is of type $(1,0)$ and the multiplication in $G$ is holomorphic, we get
$f^{-1}.df \circ I_M = i f^{-1}.df$ and hence $df \circ I_M = I_G \circ df$ for the complex structures
$I_M,I_G$ on $M$ and $G.$ This completes the proof.
\end{proof}
\subsection{Non-Abelian Cohomology}

We introduce the basic notions of non-abelian cohomology, following Onishchik \cite{Onishchik}.

\begin{definition} \label{nonabcohom}
Let $\Gamma, G$ be groups such that $\Gamma$ acts on $G$ by automorphisms. Denote this action by
$(\gamma,g) \mapsto \gamma.g.$ Let $C^1(\Gamma,G)$ denote the set of maps from $\Gamma$ to $G.$
We denote by $Z^1(\Gamma,G)$ the subset of $C^1(\Gamma,G)$ consisting of all maps $f$ satisfying
\[
f(\alpha \beta) = f(\alpha) \alpha.f(\beta) , \quad  \forall \; \alpha, \beta \in \Gamma.
\]
Note that if $\Gamma$ acts trivially, then $Z^1(\Gamma,G)= \Hom(\Gamma,G).$
The group $G$ acts on $Z^1(\Gamma,G)$ via $(g.f)(\gamma):= g f(\gamma)(\gamma.g)^{-1}.$ The set of 
$G$-orbits in $Z^1(\Gamma,G)$ is denoted by $H^1(\Gamma,G).$
If the action of $\Gamma$ is trivial, we have
\[
H^1(\Gamma,G) = \Hom(\Gamma,G) / G
\]
where $G$ acts via $(g.f)(\gamma)=g f(\gamma)g^{-1}.$
\end{definition}

\begin{definition} \label{def:actions}
Let $M$ be a smooth manifold, $G$ a Lie group, $\omega$ a $\g$-valued 1-form on $M$ and $f \in C^\infty(M,G).$ We define
\begin{equation} \label{action1}
\omega * f := \delta(f) + \Ad(f)^{-1}.\omega.
\end{equation}
If $M$ is a complex manifold and $G$ a complex Lie group, the left logarithmic derivative of  $f \in C^\infty(M,G)$ can be uniquely written as $ f^{-1}.df = f^{-1}.\partial f + f^{-1}. \bar \partial f$ with $f^{-1}.\partial f \in \Omega^{(1,0)}(M,\g)$ and $f^{-1}.\bar \partial f \in \Omega^{(0,1)}(M,\g).$ Let us now define
\begin{equation} \label{action2}
\omega \bullet f := f^{-1}.\bar \partial f + \Ad(f)^{-1}.\omega.
\end{equation}
\end{definition}

\begin{lemma} \label{curr:action}
Equation (\ref{action1}) defines a right action of $C^\infty(M,G)$ on $\Omega^1(M,\g)$ and, if $M$ and $G$ are complex, an action of $\Hol(M,G)$ on $\Omega^1_h(M,\g).$ If $\omega$ satisfies the MC equation, then so does~$\omega * f.$
Similarly, equation (\ref{action2}) defines a right action of $C^\infty(M,G)$ on $\Omega^{(0,1)}(M,\g)$ which leaves the set of solutions of the HMC equation invariant.
\end{lemma}
\begin{proof}
See Appendix \ref{App:currAction}.
\end{proof}

\begin{definition}
Let $\ZDR(M,\g)$ and $\ZDRhol(M,\g)$ denote the $\g$-valued 1-forms, resp. the holomorphic $\g$-valued 1-forms
on $M$ satisfying the MC equation. We write
\[
\HDR(M,\g) := \ZDR(M,\g) / C^\infty(M,G)
\]
and 
\[
\HDRhol(M,\g) := \ZDRhol(M,\g) / \Hol(M,G)
\]
for the sets of orbits of action (\ref{action1}) from Definition \ref{def:actions}.
\end{definition}

\subsection{Flat Holomorphic Bundles}

In this subsection, non-abelian cohomology and the period homomorphism are combined to construct a sequence of pointed sets including the set of flat holomorphic principal $G$-bundles. The results in this subsection can be found in \cite{Onishchik}. 

\begin{definition}
Let $q: \widetilde M \rightarrow M$ be the universal cover of a smooth manifold $M,$ and let $G$ be a Lie group. Given $\phi \in Z^1(\pi_1(M), \Hol(\widetilde M, G)),$ we define 
\[
P_\phi := (\widetilde M \times G) / \pi_1(M),
\]
where $\pi_1(M)$ acts by $\gamma.(m,g) = (\gamma.m , \phi(\gamma)(\gamma.m)g).$ Let the projection $p: P_\phi \rightarrow M$ be given by $p(\left[ (m,g) \right]) = q(m).$ 
\end{definition}
\begin{lemma} \label{homBdle}
The bundle $p: P_\phi \rightarrow M$ is a holomorphic principal $G$-bundle. The pullback of $P_\phi$ under $q$ is trivial. Conversely, if $p: P \rightarrow M$ is any holomorphic principal $G$-bundle such that its pullback under $q$ is trivial, then the bundle is equivalent to some $P_\phi.$ We get equivalent bundles $P_\phi \cong P_\psi$ if and only if $\phi$ and $\psi$ lie 
in the same orbit under the action of $\pi_1(M),$ i.e. if and only if $\left[ \phi \right]=\left[ \psi \right]$
in $H^1(\pi_1(M), \Hol(\widetilde M, G)).$
\end{lemma}

\begin{lemma} \label{lem:classBun}
Let $P$ be the map that assigns to $\omega \in \ZDRhol(M,\g)$ its period homomorphism. Let $\flhoBun(M,G)$ denote the set of flat holomorphic principal $G$-bundles over $M,$ modulo holomorphic equivalence. In this set, we choose the class of the trivial bundle as base point. Then we get an exact sequence of pointed sets
\[
G \rightarrow \Hol(M,G) \stackrel{\delta}{\rightarrow} \ZDRhol(M,\g) \stackrel{P}{\rightarrow} 
\Hom(\pi_1(M),G) \rightarrow \flhoBun(M,G) \rightarrow 0 
\]
which is exact in the sense that the fibers of $P$ coincide with the orbits of $\Hol(M,G)$ in $\ZDRhol(M,G).$
Exactness at $\Hom(\pi_1(M),G)$ means that the bundle corresponding to $\phi$ is holomorphically trivial if and only if 
$\phi = \mathrm{per}_\omega$ for some $\omega \in \ZDRhol(M,\g).$ The sequence induces another exact sequence:
\[
H^1_{\mathrm{dR,h}}(M,\g) \stackrel{P}{\rightarrow} H^1(\pi_1(M),G) \rightarrow \flhoBun(M,G) \rightarrow 0 
\]
\end{lemma}

\begin{lemma}
The map $P$ in the above sequence is surjective if one of the following conditions holds:
\begin{itemize}
\item $\pi_1(M)$ is a free group and $G$ is connected.
\item $\pi_1(M)$ is a free abelian group and $G$ is a compact connected group whose
			cohomology is torsion free.
\end{itemize}
In particular, under these conditions every flat holomorphic principal $G$-bundle is holomorphically trivial.
\end{lemma}

\subsection{Trivial Bundles}

If $P = M \times G$ is a smoothly trivial principal $G$-bundle, then the gauge group of $P$ can be identified with the group $C^\infty(M,G)$ of $G$-valued smooth maps on $M$ under pointwise multiplication. This will be made precise in the following lemma.

\begin{lemma} \label{lem:autPbracket}
If $P=M \times G \rightarrow M$ is a trivial principal $G$-bundle, then we have semidirect product decompositions
\[
\Aut(P) = \diff(M) \ltimes C^\infty(M,G) \quad \text{and} \quad \aut(P) = \vect(M) \ltimes C^\infty(M, \g),
\]
and the Lie bracket of $\vect(M) \ltimes C^\infty(M, \g)$ is given by
\begin{equation} \label{eq:autPbracket}
[(X,f),(Y,g)]= (-[X,Y], Y(f) -X(g) + [f,g]).
\end{equation}
\end{lemma}
\begin{proof}(Sketch)
If $\phi(x,g)=(\phi_1(x,g),\phi_2(x,g))$ is an automorphism of $P,$ then $\phi(x,g)=\phi(x,1)g=(\phi_1(x,1),\phi_2(x,1))g,$ and it is easy to see that $\gamma(x) = \phi_1(x,1) \in \diff(M)$ and $f(x)= \phi_2(x,1) \in C^\infty(M,G).$ Now given $\phi_i=(\gamma_i, f_i),$ one computes that
\[
\phi_1 \circ \phi_2(x,g)= (\gamma_1 \circ \gamma_2(x), f_1(\gamma_2(x))f_2(x)g).
\]
Hence $\Aut(P) = \diff(M) _\alpha \ltimes C^\infty(M,G),$ where the action $\alpha$ of $\diff(M)$ on $C^\infty(M,G)$ is given by $\phi.f(x)= f(\phi^{-1}(x)).$ The corresponding derived action of $\vect(M)$ on $C^\infty(M, \g)$ is then computed to be $X.f = -X(f),$ and equation (\ref{eq:autPbracket}) follows from the general formula for the Lie bracket of a semidirect product.
\end{proof}

\begin{remark}
The Lie bracket of two vector fields $X$ and $Y,$ if we think of $\vect(M)$ as Lie algebra of $\diff(M),$ is given by $-[X,Y].$ This is the reason for the occurence of $-[X,Y]$ in equation (\ref{eq:autPbracket}).
\end{remark}

\subsection{Bundles over Riemann Surfaces}

\begin{definition}
From now on, $\Sigma$ will denote a \emph{closed Riemann surface}, i.e. a compact complex manifold of complex dimension 1, and let $g$ denote the \emph{genus} of $\Sigma.$ If $H_\Sigma = \Omega^1_h(\Sigma)$ denotes the vector space of holomorphic 1-forms on $\Sigma,$ then $g = \dim_\C H_\Sigma.$   
\end{definition}

\begin{theorem} \label{thm:topBdles}
Let $\Sigma$ be a closed Riemann surface of genus $g$ and $G$ a connected topological group.
Then there is a bijective correspondence between the set of isomorphism classes of topological principal 
$G$-bundles on $\Sigma$ and $\pi_1(G),$ the fundamental group of $G.$
\end{theorem}
\begin{proof}
See \cite{Ramanathan:1975}, Proposition 5.1.
\end{proof}
Principal $G$-bundles are determined by homotopy classes of maps into the classifying space $BG$ of $G,$ and since each homotopy class has a smooth representative (see the discussion in Chapter 4 of \cite{Hirsch:1994}), we have that a topologically trivial bundle is also smoothly trivial. Thus we get a corollary of Theorem \ref{thm:topBdles}:

\begin{corollary} \label{cor:Gsimpconn}
If $G$ is simply connected, then every principal $G$-bundle over $\Sigma$ is smoothly trivial.
\end{corollary}
\subsubsection{Riemann Surfaces of Genus 1}

If $g=1,$ then the fundamental group $\pi_1(\Sigma)$ is isomorphic to $\Z^2.$ This implies that a homomorphism $\pi_1(\Sigma) \rightarrow G$ is uniquely determined by an ordered pair of commuting elements of $G.$ 

\begin{theorem}
Let $K$ be a compact simply connected Lie group and $\sigma$ an automorphism of $K.$
Then the set $F$ of fixed points of $\sigma$ is connected.
\end{theorem}
\begin{proof}
See \cite{Borel:1961}, Theorem 3.4.
\end{proof}

\begin{corollary}
The centralizer $Z(g)$ of any element $g \in K$ is connected.
\end{corollary}
\begin{proof}
$Z(g)$ is the set of fixed points of the automorphism $I_g: x \mapsto gxg^{-1}.$
\end{proof}

\begin{corollary}
Let $K$ be a compact simply connected Lie group and $x$ and $y$ two commuting elements of $K.$ Then there exists a maximal torus $T$ in $K$ containing both $x$ and $y.$
\end{corollary}
\begin{proof}
By the previous corollary, $Z(x)$ is a connected subgroup of $K.$ It is closed, thus
compact, and it contains $y$ by hypothesis. So $y$ is contained in a maximal torus
$T'$ of $Z(x),$ which is also a torus in $K.$ Every torus $T' \subset K$ is contained in some maximal torus $T$
of $K,$ and this maximal torus $T$ contains both $x$ and $y.$
\end{proof}

\begin{corollary}
Let $K$ be a compact simply connected Lie group, $T$ a maximal torus in $K,$ and $W$ the Weyl group of $K.$ Then the space of holomorphic principal $K$-bundles over a closed Riemann surface of genus 1 is isomorphic to 
$\Hom(\Z^2, T)/W,$ where $W$ acts on $T$ in the usual way. As a set, $\Hom(\Z^2, T)$ is isomorphic to $T \times T.$
\end{corollary}

\begin{remark}
This corollary is Theorem 2.6 in the article \cite{FriedmanMorganWitten}, in which Friedman, Morgan and Witten investigate principal $G$-bundles over elliptic curves. They describe the set $\Hom(\Z^2, T)/W$ as a weighted projective space, the weights being determined by the Weyl group $W.$ See also \cite{EtingofFrenkel}, Proposition 4.1.
\end{remark}
%
\section{Complex Structures on Principal Bundles}

In this section we describe almost complex structures on the total space $P$ of a smoothly trivial principal $G$-bundle with $G$ being a complex Lie group and the base manifold $M$ a complex manifold. We only consider almost complex structures for which there is a bundle atlas consisting of holomorphic maps, and call these \emph{compatible almost complex structures.} Thus, if a compatible structure on $P$ is integrable, then $P$ is a holomorphic principal bundle. We describe the integrability condition and then we state the condition under which two such holomorphic bundle structures are isomorphic. An analogous discussion of holomorphic structures on \emph{vector} bundles can be found for example in \cite{Kobayashi}. A common feature of the vector bundle and the principal bundle case is that the space of almost complex structures is an affine space modeled on a vector space of forms of type $(0,1).$ In case of a vector bundle $E,$ the forms take values in the bundle $\End(E),$ whereas in case of a principal $G$-bundle, they are $\g$-valued forms.

Let $P$ be a smoothly trivial holomorphic principal $G$-bundle. Then $P$ is a complex manifold, and it carries the corresponding natural almost complex structure.Since the bundle charts are biholomorphic, they can be used to express this almost complex structure in terms of a certain $\g$-valued form on $M.$

 \begin{lemma} \label{lem:locStruct}
Suppose $P$ is a smoothly trivial holomorphic principal $G$-bundle over $M,$ and $\psi_j:P|_{U_j} \rightarrow U_j \times G$ a holomorphic bundle chart given by $\psi_j(x,g)=(x, \gamma_j(x)g)$ for smooth maps $\gamma_j: U_j \rightarrow G.$ Since $P$ is a complex manifold, it carries an almost complex structure $I,$ and on $P|{U_j}$ this structure is given by
\begin{equation} \label{eq:locStruct}
I_{(x,g)}(v,w) = (iv, iw + 2i \gamma_j(x)^{-1}.\bar \partial_x \gamma_j(v).g).
\end{equation}
There is a form $\omega \in \Omega^{(0,1)}(M,\g)$ which is locally given by $\gamma_j^{-1}. \bar \partial \gamma_j.$
\end{lemma}
\begin{proof}
See Appendix \ref{AppA:struct}.
\end{proof}

Conversely, if we are given a $\g$-valued form $\omega$ of type $(0,1),$ an elementary computation shows that equation (\ref{eq:locStruct}) defines an almost complex structure on the product $M \times G.$

\begin{definition}
If $\omega \in \Omega^{(0,1)}(M, \g),$ let $I_\omega$ be the almost complex structure on $M \times G$ defined by $I_\omega(v,w) = (iv, iw + 2i \omega_x(v).g)$ for $(v,w) \in T_{(x,g)}(M \times G).$
\end{definition}

The equation $\omega = \gamma^{-1}.\bar \partial \gamma$ is locally solvable if and only if $\omega$ satisfies the HMC equation (see \cite{Onishchik}, Lemma 6.1). If that is the case, then we can cover $M$ with open sets $U_i$ on which there is a solution $\gamma_i,$ and then the $\gamma_i$ can be used to define a holomorphic bundle atlas on $M \times G.$ Recall the characterization of integrability in terms of torsion in Theorem \ref{integrability}.  We will use this characterization to prove that $I_\omega$ is integrable if and only if $\omega$ satisfies the HMC equation.

\begin{theorem} \label{MainInt}
Let $M$ be a complex manifold, $G$ a complex Lie group and let $\omega \in \Omega^{(0,1)}(M, \g).$ The compatible almost complex structure $I_\omega$ on the smoothly trivial principal $G$-bundle $M \times G$ is integrable if and only if $\omega$ satisfies the HMC equation.
\end{theorem}
\begin{proof}
See Appendix \ref{AppA:Int}.
\end{proof}

This theorem shows that to each solution $\omega$ of the HMC equation there is a smoothly trivial holomorphic principal $G$-bundle over $M.$ Now we want to show that two such bundles are holomorphically equivalent if the corresponding forms are in the same orbit under the action (\ref{action2}) in Definition \ref{def:actions}.

\begin{theorem} \label{thm:equiv}
Let $\omega_1$ and $\omega_2$ be solutions to the HMC equation, and let $P_1$ and $P_2$ be the corresponding smoothly trivial holomorphic principal $G$-bundles over $M.$ Then $P_1$ and $P_2$ are holomorphically equivalent if and only if $\omega_1 = \omega_2 \bullet f$ for some $f \in C^\infty(M,G),$ i.e. if and only if $\omega_1$ and $\omega_2$ are in the same orbit under the action (\ref{action2}) of $C^\infty(M,G).$
\end{theorem}
\begin{proof}
See Appendix \ref{AppA:Equiv}.
\end{proof}

\begin{remark}
If $M= \Sigma$ is a closed Riemann surface, then $\Omega^{(0,2)}(\Sigma) = 0,$ hence every form of type $(0,1)$ satisfies the HMC equation. If furthermore $G$ is simply connected, then by Corollary \ref{cor:Gsimpconn} all principal $G$-bundles are smoothly trivial, and hence the set of isomorphism classes of  holomorphic principal $G$-bundles is given by $\Omega^{(0,1)}(\Sigma, \g)/C^\infty(\Sigma,G).$
\end{remark}
%
\section{Central Extensions and Coadjoint Orbits}

In this final section we want to establish a connection between holomorphic principal $G$-bundles over a closed Riemann surface $\Sigma$ and coadjoint orbits in the dual of a central extension of the Lie algebra  $C^\infty(\Sigma, \g).$ The key ingredient is a pairing which involves $\g$-valued forms of type $(0,1)$ on the one hand, and the algebra $C^\infty(\Sigma, \g)$ on the other hand. The pairing is defined as integral of an appropriate 2-form over $\Sigma.$ This material was published in \cite{EtingofFrenkel} and generalizes a similar construction for loop groups (see \cite{PressleySegal}, Chapter 4.3).

\subsection{Central Extensions of the Current Algebra}
We recall some basic facts about extensions of Lie algebras. 
\begin{definition}
A short exact sequence $ 0 \rightarrow \mathfrak i \stackrel{\iota}{\rightarrow} \mathfrak g \stackrel{\pi}{\rightarrow} \mathfrak h \rightarrow 0$ of Lie algebras is called \emph{extension of $\mathfrak h$ by $\mathfrak i$.} Sometimes we just say that $\g$ is the extension, if $\iota$ and $\pi$ are understood. The extension is called \emph{central} if $\mathfrak i$ is abelian and is mapped by $\iota$ into the center of $\g.$ Given a second extension $ 0 \rightarrow \mathfrak i \stackrel{\iota'}{\rightarrow} \mathfrak g' \stackrel{\pi'}{\rightarrow} \mathfrak h \rightarrow 0$ of $\mathfrak h$ by $\mathfrak i,$ a \emph{morphism} or \emph{isomorphism} of extensions is given by a  morphism or isomorphism $\phi: \g \rightarrow \g'$, respectively, such that $\phi \circ \iota = \iota'$ and $\pi = \pi' \circ \phi.$
\end{definition}

\begin{definition}
Let $Z^2(\mathfrak h, \mathfrak i)$ denote the vector space of bilinear maps $\Omega: \mathfrak h \times \mathfrak h \rightarrow \mathfrak i$ such that $\Omega(v,v) = 0$ and $\Omega([u,v],w)+ \Omega([v,w],u)+ \Omega([w,u],v) = 0.$ These are the Lie-algebra two-cocycles. For a linear map $\beta: \mathfrak h \rightarrow \mathfrak i,$ let $\delta(\beta)(u,v)= -\beta([u,v]).$ Then $\delta(\beta) \in Z^2(\mathfrak h, \mathfrak i),$ and the image under $\delta$ of all linear maps is called $B^2(\mathfrak h, \mathfrak i),$ the space of two-coboundaries.
\end{definition}

\begin{lemma}
If $\Omega \in Z^2,$ we define a central extension of $\mathfrak h$ by $\mathfrak i$ as $\g = \mathfrak i \oplus \mathfrak h$ with Lie bracket $[(u,x),(v,y)]= (\Omega(x,y),[x,y]).$ Two such central extensions given by $\Omega$ and $\Omega'$ are equivalent if and only if $\Omega -\Omega'$ is a coboundary. The vector space of equivalence classes of central extensions of $\mathfrak h$ by $\mathfrak i$ is isomorphic to the second Lie algebra cohomology $H^2(\mathfrak h, \mathfrak i)= Z^2(\mathfrak h, \mathfrak i)/B^2(\mathfrak h, \mathfrak i)$ of $\mathfrak h$ with values in $\mathfrak i.$
\end{lemma}

Now let us consider Lie algebra $C^\infty(M, \g).$ Assume that $\g$ is simple and let $\left\langle \cdot, \cdot \right\rangle$ be an invariant inner product on $\g.$ Then if $f, g \in C^\infty(M,\g),$ we have $dg \in \Omega^1(M,\g)$ and we define the $\g$-valued 1-form $\left\langle f, dg \right\rangle$ by forming the inner product pointwise. The following result is Proposition 4.2.8. in \cite{PressleySegal}.
\begin{lemma}
The assignment
\[
(f,g) \mapsto \left\langle f, dg \right\rangle
\]
defines a cocycle of $C^\infty(M,\g)$ with values in $\Omega^1(M,\g)/dC^\infty(M,\g).$ The corresponding central extension is \emph{universal} in the sense that it has every other central extension as a quotient. The extensions of $C^\infty(M, \g)$ by $\R$ correspond to one-dimensional closed currents $C$ on $M,$ the cocycle being given by integrating above cocycle over $C.$
\end{lemma}

\begin{definition}
Let $\Sigma$ be a closed Riemann surface of genus $g,$ and let $G$ be a simple Lie group and $\eta \in H_\Sigma.$
Let $E_\eta$ be the central extension of $C^\infty(\Sigma, \g)$ by $\R$ given by the cocycle
\[
\Omega_\eta(f,g) = \int_\Sigma \eta \wedge \langle f, dg \rangle.
\]
\end{definition}

\begin{lemma} \label{lem:ExtAdj}
The group $C^\infty(\Sigma, G)$ acts on $E_\eta$ by
\[
f.(x,g)= \left( x- \int_\Sigma \eta \wedge  \langle f^{-1}. \bar  \partial f, g \rangle,  \, \Ad(f).g  \right)
\]
\end{lemma}
\begin{proof}
See \cite{MaierNeeb}, Proposition III.3.
\end{proof}
%
%
%

\subsection{Coadjoint Orbits}


\begin{definition} \label{smoothDual}
Let  $E^* = \C \oplus \Omega^{(0,1)}(\Sigma,\g),$ and define a right action of $C^\infty(\Sigma, G)$ on $E^*$ by
\begin{equation} \label{eq:orbAction}  
(\lambda ,\xi)*f = (\lambda, \lambda f^{-1}.\bar \partial f +\Ad(f)^{-1}.\xi).
\end{equation}
We define a pairing $E^* \times E_\eta \rightarrow \C$ via
\begin{equation} \label{eq:pairing}
((\lambda, \xi),(x,X)) = \lambda x - \int_{\Sigma} \eta \wedge \langle \xi,X \rangle. 
\end{equation}
\end{definition}


\begin{lemma} \label{lem:dual}
Equation (\ref{eq:orbAction}) defines a right action, and the pairing (\ref{eq:pairing}) 
is $C^\infty(\Sigma, G)$-invariant with respect to this action on $E^*$ and the action from Lemma \ref{lem:ExtAdj} 
on the extension $E_\eta.$ Furthermore, the pairing is non-degenerate.
\end{lemma}
\begin{proof}
See Appendix \ref{AppB:Pairing}.
\end{proof}

Lemma \ref{lem:dual} allows us to view $E^*$ as being a subspace of the dual of the extended Lie algebra $E_\eta,$ and the action (\ref{eq:orbAction}) as coadjoint action. Etingof and Frenkel \cite{EtingofFrenkel} call $E^*$ the \emph{smooth part} of the dual of $E_\eta.$ We conclude with the following theorem, which establishes a relation between classes of holomorphic principal bundles and orbits of the coadjoint action.

\begin{theorem} \label{MainThm}
There is a one-to-one correspondence between orbits of the action (\ref{eq:orbAction}) for each $\lambda \not= 0$
and equivalence classes of holomorphic principal $G$-bundles over $\Sigma$.
\end{theorem}
\begin{proof}
See Appendix \ref{AppB:Corresp}.
\end{proof}
%
\appendix
\section{Appendix: Technical Details}

%
\subsection{Gauge Group Actions} \label{App:currAction}

Let $\omega$ be a smooth (holomorphic) 1-form and $f,g$ smooth (holomorphic) $G$-valued maps. The product rule $\delta(fg) = \Ad(g)^{-1}.\delta(f) + \delta(g)$ is easily verified. We use it to get
\begin{equation*} \begin{split}
\omega * (fg) &= \delta(fg) + \Ad(fg)^{-1}.\omega                                        \\
							&= \Ad(g)^{-1}.\delta(f) + \delta(g) + \Ad(g)^{-1} \Ad(f)^{-1}. \omega      \\
							&= \delta(g) + \Ad(g)^{-1}.(\delta(f)+ \Ad(f)^{-1}.\omega)                  \\
							&= (\omega*f)*g.
\end{split} \end{equation*} 
Because for constant functions $f$ the left logarithmic derivative $\delta(f)$ is zero, the 
identity element of the two mapping groups acts trivially. It remains to show that if $M$ and $G$ are 
complex and $f$ and $\omega$ are holomorphic, then $\omega*f$ is a holomorphic 1-form.
First observe that, because on $G$ multiplication and inversion are holomorphic, we get that 
$f^{-1}.df \circ I = if^{-1}.df$ and similarly for $\Ad(f)^{-1}. \omega,$ so that $\omega *f$ is of
type $(1,0).$ We will need the following equation:
\begin{equation} \label{ActMC}
d(\omega *f)+ \frac{1}{2} \left[ \omega *f,\omega *f \right] 
= \Ad(f)^{-1}.\left(d\omega+\frac{1}{2} \left[ \omega ,\omega \right] \right)
\end{equation}
\begin{proof}(of Equation (\ref{ActMC})).
We use equation
\[
d(\Ad(f)^{-1}. \omega) = \Ad(f)^{-1}.d \omega - [ \delta(f), \Ad(f)^{-1}. \omega ] 
\]
from the proof of \cite{MaierNeeb}, Proposition III.3. (see also \cite{Sharpe} Chapter 3, Exercise 4.14) and get
\begin{equation*} \begin{split}
d(\omega *f)+ \frac{1}{2} \left[ \omega *f,\omega *f \right] &= d(\delta(f)) + \Ad(f)^{-1}.d\omega - [\delta(f), \Ad(f)^{-1}. \omega] + \frac{1}{2}[\delta(f), \delta(f)]  \\
& \quad + [\delta(f), \Ad(f)^{-1}. \omega] + \frac{1}{2}[\Ad(f)^{-1}.\omega, \Ad(f)^{-1}.\omega] \\
&= \Ad(f)^{-1}.(d \omega + \frac{1}{2}[\omega, \omega]).
\end{split} \end{equation*}
\end{proof}
The right hand side of (\ref{ActMC}) is of type $(2,0)$ because of Remark \ref{rem:typeTwozero} and because $\omega$ and $f$  are holomorphic. The left hand side has to be of type $(2,0)$ too, and this is only possible if $\bar \partial (\omega * f)=0.$ Hence $\omega *f$ is holomorphic.

It is immediate from Equation (\ref{ActMC}) that $\omega *f$ satisfies the MC equation if and only if $\omega$ does.

Now let us consider the action denoted by $\omega \bullet f.$ The equality $\omega \bullet (fg) = (\omega \bullet f) \bullet g,$ as well as the fact that $\omega \bullet f$ is again of type $(0,1),$ follow as above.
We need to show that the action leaves the set of solutions of the HMC equation invariant. So assume that 
$\bar \partial \omega + \frac{1}{2}\left[ \omega, \omega \right] =0.$ Observe that $\omega * f = \omega \bullet f + f^{-1}.\partial f.$ If we plug this into equation (\ref{ActMC}) and expand both sides, the left hand side yields
\begin{multline*}
\partial(\omega \bullet f) + \partial(f^{-1}.\partial f) + \bar \partial (\omega \bullet f) 
+ \bar \partial (f^{-1}.\partial f)   \\
+ \left[\omega \bullet f,f^{-1}.\partial f \right] + \frac{1}{2}\left[\omega \bullet f,\omega \bullet f \right] 
+ \frac{1}{2}\left[f^{-1}.\partial f, f^{-1}.\partial f \right]
\end{multline*}
and the right hand side yields
\[
\Ad(f)^{-1}.(\bar \partial \omega + \frac{1}{2}\left[ \omega,\omega \right] + \partial \omega).
\]
As the terms of degree (0,2) have to coincide, we get
\[
\bar \partial(\omega \bullet f) + \frac{1}{2}\left[ \omega \bullet f, \omega \bullet f \right] = 
									\Ad(f)^{-1}.(\bar \partial \omega + \frac{1}{2}\left[ \omega,\omega \right])
\]
Hence $\omega \bullet f$ satisfies the HMC equation.

%
\subsection{Compatible Almost Complex Structures} \label{AppA:struct}

We consider a bundle chart 
\[
\psi: U \times G \rightarrow U \times G , \qquad (x,g) \mapsto (x,\gamma(x)g)
\]
where $\gamma : U \rightarrow G$ is smooth.

\begin{lemma}
The differential of the bundle chart $\psi$ is given by
\[
d_{(x,g)} \psi = 
\begin{pmatrix}
\id & 0 \\
d_{\gamma(x)}\rho_g \circ d_x \gamma & d_g \lambda_{\gamma(x)}
\end{pmatrix},
\]
where $\rho_g$ and $\lambda_g$ denote multiplication by $g$ on the right and left, respectively.
\end{lemma}
\begin{proof}
For fixed $g \in G$ we let $\psi_g: x \mapsto \psi(x,g),$ and similarly we define $\psi_x: g \mapsto \psi(x,g).$ Then $\psi_g=(\id, \rho_g \circ \gamma(x))$. The chain rule yields $d_{(x,g)} \psi_g= (\id, d_{\gamma(x)} \rho_g \circ d_x \gamma).$ The differential of $\psi_x: g \mapsto (x,\lambda_{\gamma(x)}(g) )$ is  given by $d_{(x,g)} \psi_x = (0, d_g \lambda_{\gamma(x)}).$
\end{proof}

If $(M,I)$ and $(N,J)$ are almost complex manifold and $\phi: M \rightarrow N$ a biholomorphic map, then $I_xv = (d_x \phi)^{-1} (J_x d_x\phi (v)).$ We now apply this fact to the biholomorphic map $\psi: P_U \rightarrow U \times G$ where $U \times G$ carries the product almost complex structure $J(v,w)= (iv, iw).$ 

\begin{proof}(of Lemma \ref{lem:locStruct}).
We have 
\[
I_{(x,g)} = ( d_{(x,g)}\psi)^{-1} \circ J_{\psi(x,g)}  
\circ d_{(x,g)}\psi
\]
and if we let $I_{(x,g)}(v,w) = (c, d),$ then by definition of $I_{(x,g)}$ we get
\[
d_{(x,g)}\psi (c,d) = J_{(x,\gamma(x).g)} d_{(x,g)} \psi (v,w)   
\]
and
\[
(c, d_{\gamma(x)} d_x \gamma (c).g + \gamma(x).d)
 = (iv, i(d_{(d_x \gamma(a)}.g) 
+ i(\gamma(x).w)).
\]

So we see that $c =iv.$ We use this to compute $d.$ We also use the fact that
multiplication on $G$ is holomorphic, thus its differentials commute with $i.$
We get
\[
\gamma(x).d = i \gamma(x). w + i d_x \gamma (v) . g - d_x \gamma (iv) .g
\]
\[
d= iw +(i \gamma(x)^{-1} d_x \gamma(v)- \gamma(x)^{-1} d_x \gamma(iv)).g.
\]
We use $d_x\gamma(iv) = \partial_x \gamma(iv) + \bar \partial_x \gamma(iv) = i (\partial_x \gamma(v) - \bar \partial_x \gamma(v))$ and finally get
\[
(c,d)=(iv,iw + 2i \gamma(x)^{-1} \bar \partial_x \gamma(v).g),
\]
where $\omega_x(a) = \gamma(x)^{-1} \bar \partial_x \gamma(a)$ is a $\g$-valued
(0,1)-form on $U.$ 
\end{proof}


\subsection{Integrability} \label{AppA:Int}
First we recall some facts about differential forms. If $\omega$ is a 1-form and $X$ and $Y$ are vector fields on a smooth manifold $M,$ then we denote by $X\omega(Y)$ the result of applying the derivation $X$ to the smooth function $\omega(Y).$ The following holds (see \cite{Sharpe}, Chapter 1, Lemma 5.15):
\begin{equation*}
d\omega(X,Y)= X\omega(Y)-Y\omega(X) - \omega([X,Y]).
\end{equation*}
Now if $M$ is a complex manifold and $\omega$ is of type $(0,1),$ then we also have
\begin{equation*}
d\omega(X,Y) - d\omega(iX,iY) = 2 \bar \partial \omega(X,Y).
\end{equation*}
We can combine these two equations to get
\begin{equation*} \begin{split} 
2 \bar \partial \omega(X,Y) &= X\omega(Y)-Y\omega(X) - \omega([X,Y]) -(iX\omega(iY)-iY\omega(iX) - \omega([iX,iY]))  \\
&= X\omega(Y)-Y\omega(X) - \omega([X,Y]) -(X\omega(Y)-Y\omega(X) + \omega([X,Y]))  \\
&= -2\omega([X,Y]),
\end{split} \end{equation*}
and we divide by 2 to get
\begin{equation} \label{consequence}
\bar \partial \omega(X,Y)=-\omega([X,Y]).
\end{equation}

The following lemma constitutes the main part in our proof of Theorem \ref{MainInt}
\begin{lemma} \label{lem:Torsion}
Equip the manifold $M \times G$ with the almost complex structure $I_\omega$ corresponding to $\omega \in \Omega^{(0,1)}(M, \g),$ and let $N_\omega$ be the torsion of $I_\omega.$ Identify $\aut(M \times G) \cong \vect(M) \ltimes C^{\infty}(M,\g)$ as in Lemma \ref{lem:autPbracket}. Then for $X,Y \in \vect(M),$  the torsion of $(X,0)$ and $ (Y,0)$ is given by 
\begin{equation*} \begin{split}
N_\omega((X,0),(Y,0)) &= (0,-2\left[ \omega,\omega \right](X,Y)-4 \bar \partial \omega(X,Y)) \\
&= (0, -4 \bar F(\omega)(X,Y)), 
\end{split} \end{equation*}
which vanishes for all $X,Y$ if and only if $\omega$ satisfies the HMC equation.
\end{lemma}
\begin{proof}
By definition of $I_\omega$ and by the formula for the Lie bracket from Lemma \ref{lem:autPbracket}, the first component of the torsion will simply be the torsion of the natural almost complex structure on $M,$ which vanishes since $M$ is a complex manifold. Also, the second term $[(X,0),(Y,0)]$ has zero as second component and thus plays no role in our computation. Let us compute the first, third and fourth term of the torsion.
\\
{\bf First term: $[I_\omega(X,0),I_\omega(Y,0)]$}
\begin{equation*}
I_\omega(X,0) = (iX, 2i\omega(X)).
\end{equation*} 
Now we compute the Lie bracket
\begin{equation} \begin{split} \label{FirstTerm}
[(iX, 2i\omega(X)),(iY, 2i\omega(Y))]  &= (-[iX, iY], iY(2i\omega(X))-iX(2i\omega(Y))+4[i\omega(X), i\omega(Y)] ) \\
&= (-[iX, iY], -2iY(\omega(iX))+2iX(\omega(iY))+4[\omega(iX), \omega(iY)] )
\end{split} \end{equation} 

{\bf Third term: $I_\omega[I_\omega(X,0),(Y,0)]$}
As above, $I_\omega(X,0)=(iX, 2i\omega(X))$ and then
\begin{equation*}
[(iX, 2i\omega(X)),(Y,0)]= (-[iX, Y], Y(2i\omega(X))).
\end{equation*} 
Now we need to apply $I_\omega.$
\begin{equation} \begin{split} \label{ThirdTerm}
I_\omega(-[iX, Y], Y(2i\omega(X))) &= (-i[iX,Y], iY(2i\omega(X))-2i\omega([iX,Y])) \\
&= (-i[iX,Y], -2iY(\omega(iX))+2\omega(i[iX,Y]))
\end{split} \end{equation} 

{\bf Fourth term: $I_\omega[(X,0),I_\omega(Y,0)]$}
Clearly, we just need to switch $X$ and $Y$ in Equation \ref{ThirdTerm} and  multiply by $-1.$ Thus we get
\begin{equation} \label{FourthTerm}
2iX(\omega(iY))-2\omega(i[iY,X])
\end{equation} 
as second component.

Finally  we compute the second component of the torsion $N_\omega((X,0),(Y,0))$ by subtracting (\ref{ThirdTerm}) and (\ref{FourthTerm}) from (\ref{FirstTerm}) to get
\begin{gather*}
-2iY(\omega(iX))+2iX(\omega(iY))+ 4[\omega(iX), \omega(iY)] +2iY(\omega(iX))- 2\omega(i[iX,Y]) \\
\qquad- 2iX(\omega(iY))+ 2\omega(i[iY,X])
\end{gather*}
We have some cancellation, and then we use $\C$-bilinearity of the Lie bracket on $\g$ to get
\begin{equation*} \begin{split}
-2[\omega, \omega](X,Y) +2\omega([X,Y]) - 2\omega([Y,X]) = -2[\omega, \omega](X,Y) +4 \omega([X,Y]).
\end{split} \end{equation*} 
Now by Equation (\ref{consequence}), the latter term equals $-2[\omega, \omega](X,Y) -4 \bar \partial \omega(X,Y)= -4 \bar F(\omega)(X,Y).$
\end{proof}

This lemma shows that $\bar F(\omega) =0 $ is a necessary condition for the structure $I_\omega$ to be integrable. To see that it is also sufficient, note that the torsion is $C^\infty(M \times G)$-bilinear and symmetric (see \cite{KobayashiNomizu:1963}). Hence it suffices to compute the torsion of pairs $(X,0),(Y,0)$ and $(X,0),(0,f)$ and $(0,f),(0,g)$ of elements of $\vect(M) \ltimes C^\infty(M,\g)$. By Lemma \ref{lem:Torsion}, the torsion $N((X,0),(Y,0))$ vanishes, and the two remaining expressions $N((X,0),(0,f))$ and $N((0,f),(0,g))$ can easily be seen to vanish, too. This completes the proof of Theorem \ref{MainInt}.


\subsection{Equivalence of Complex Structures} \label{AppA:Equiv}
\begin{proof}
An isomorphism of holomorphic principal $G$-bundles is a biholomorphic bundle map 
$F:M \times G \rightarrow M \times G.$ Every such bundle map is of the form $(x,g) \mapsto (x,f(x)g)$ 
for a smooth map $f \in C^\infty(M,G).$
Now there is an open cover $\{ U_j \}$ of $M$ such that we have bundle charts $\Psi_j$ and $\Psi'_j$ for the two bundles,
given by
\[
(x,g) \mapsto (x, \psi_j(x)g) ,\quad (x,g) \mapsto (x, \gamma_j(x)g)
\]
respectively, and such that the restrictions of the forms defining the structures are given by
$\xi_1|U_j=\psi_j^{-1} \bar \partial \psi_j$ and $\xi_2|U_j=\gamma_j^{-1} \bar \partial \gamma_j.$
Let us omit the indices as we will stay within one chart $U=U_j.$ Then the map $F$ is given in this chart by $  \Psi' \circ F \circ \Psi^{-1}: (x,g) \mapsto (x,\gamma(x) f(x) \psi(x)^{-1}g).$
This map has to be holomorphic with respect to the product complex structure on $U \times G,$
and thus $x \mapsto \gamma(x) f(x) \psi(x)^{-1}$ has to be holomorphic. Using the product rule for $G$-valued functions, we get
\begin{eqnarray*}
0 &=&  \bar \partial (\gamma \cdot f \cdot \psi^{-1})  \\  
	&=&  \bar \partial \gamma .f. \psi^{-1} + \gamma .\bar \partial f. \psi^{-1} - 
	\gamma.f.\psi^{-1}.\bar \partial \psi .\psi^{-1}.
\end{eqnarray*}      
Now multiplying by $\gamma^{-1}$ from the left and $ \psi.f^{-1}$ from the right yields
\[
0= \gamma^{-1} \bar \partial \gamma + \bar \partial f. f^{-1} - f.\psi^{-1}.\bar \partial \psi.f^{-1} 
\]
and thus
\[
\xi_1 = f^{-1}. \bar \partial f + \Ad(f)^{-1}. \xi_2,
\]
which proves one part of our claim. Conversely, let $\xi_1 = f^{-1}. \bar \partial f + \Ad(f)^{-1}. \xi_2$ such that in a coordinate patch $U$ we have $\xi_1=\psi^{-1}.\bar \partial \psi$ and  $\xi_2=\phi^{-1}.\bar \partial \phi.$ Then by reversing the above computation we get $\bar \partial (\gamma \cdot f \cdot \psi^{-1})=0,$ and therefore $\Psi:(m,g) \rightarrow (m,f(m)g)$ is a biholomorphic equivalence between the bundles corresponding to $\xi_1$ and $\xi_2.$
\end{proof}

%
 

\subsection{Action and Pairing} \label{AppB:Pairing}

\begin{proof}
By Lemma~\ref{curr:action}, equation (\ref{eq:orbAction}) defines a right action. We need to prove $C^\infty(\Sigma, G)$-invariance of the pairing, so we need to show
\[
((\lambda \bar \partial + \xi )*f, (x,X))=(\lambda \bar \partial + \xi, f.(x,X))
\]
for all $f \in C^\infty(\Sigma,G), X \in C^\infty(\Sigma, \g)$ and $\xi \in \Omega^{(0,1)}(\Sigma,\g).$
Computing the left hand side yields:
\begin{equation*} \begin{split}
((\lambda \bar \partial + \xi )*f, (x,X)) &= (\lambda \bar \partial + \lambda f^{-1}. \bar \partial f +                                                                   \Ad(f)^{-1}.\xi, (x,X)) \\
																					&= \lambda x - \int_\Sigma \eta \wedge
																							\langle \lambda f^{-1}. \bar \partial f + \Ad(f)^{-1}.\xi, X \rangle 
\end{split} \end{equation*}
and the right hand side is:
\begin{equation*} \begin{split}
(\lambda \bar \partial + \xi, f*(x,X)) &= \left( \lambda \bar \partial + \xi, \left( x -\int_\Sigma \eta \wedge 
																						\langle f^{-1}.\bar \partial f, X \rangle , \Ad(f).X \right) \right)  \\
							&=\lambda \left( x -\int_\Sigma \eta \wedge \langle f^{-1}.\bar \partial f, X \rangle \right)
																								- \int_\Sigma \eta \wedge \langle \xi, \Ad(f).X \rangle
\end{split} \end{equation*}
which coincides with the left hand side since $\langle \cdot, \cdot \rangle $ is invariant.	

Now it remains to prove that the pairing is non-degenerate. Here we will use that integration over manifolds is defined using integration over coordinate patches together with a partition of unity.  The main fact is that on an open set $U \subset \C$ the holomorphic form $\eta$ is given by $fdz$ for a holomorphic nonzero function $f.$  Similarly a form of type $(0,1)$ is locally of the form $gd \bar z$ for a smooth function $g$ on $U$ with values in $\g.$ 
Here $dz_u$ and $d\bar z_u$ are at every point $u \in U$ a basis of the one-dimensional spaces 
$T_u\Sigma_{(1,0)}$ and  $T_u\Sigma_{(0,1)}$ (recall Definition \ref{def:Decomp}). 
So the integral of $\eta \wedge \langle \xi, X \rangle$ over $U$ becomes
\[
\int_U f \langle g,X \rangle dz \wedge d\bar z.
\] 
The function $f$ is holomorphic on $U,$ and we can assume that $U$ is connected, so that $\mathrm{supp}(f|U)=U.$
So if this integral vanishes for every smooth $X: U \rightarrow \g,$ necessarily the map $g$ has to be 
zero. Conversely if the integral vanishes for every $g,$ then $X =0.$ From this it follows 
that the pairing is non-degenerate.
\end{proof}

\subsection{Orbit-Bundle Correspondence} \label{AppB:Corresp}
\begin{proof}
The hyperplanes $\lambda = const.$ are fixed under the action. So given 
$\xi \in \Omega^{(0,1)}(\Sigma,\g)$ we want to construct a principal $G$-bundle. 
Note that on $\Sigma$ there are no non-zero (0,2)-forms, therefore 
$\left[ \xi, \xi \right]+ 2\bar \partial \xi =0,$ that is, $\xi$ satisfies the HMC equation. Using Theorem \ref{MainInt}, we get the structure of a complex manifold on 
$\Sigma \times G$ with charts 
\[
\psi_i: U_i \times G \rightarrow U_i \times G , \qquad (x,g) \mapsto (x,\gamma_i(x)g) 
\] 
and the $\gamma_i$ satisfying $\gamma_i^{-1} .\bar \partial \gamma_i = \xi$ on $U_i.$ 
This implies that, on the intersection of $U_i$ and $U_j,$ we have 
$\gamma_j^{-1} .\bar \partial \gamma_j = \gamma_i^{-1} .\bar \partial \gamma_i$ and therefore 
\[ 
\bar \partial \gamma_i = \gamma_i \gamma_j^{-1} .\bar \partial \gamma_j .  
\] 
Using the quotient rule for differentials of $G$-valued functions, we see that the transition functions 
$g_{ij}=\gamma_i \gamma_j^{-1}$ are holomorphic:
\begin{equation*}  \begin{split} 
\bar \partial (\gamma_i \gamma_j^{-1})   
		&=  (\bar \partial \gamma_i )\gamma_j^{-1} - \gamma_i \gamma_j^{-1} (\bar \partial \gamma_j )\gamma_j^{-1}  \\ 
		&= \gamma_i \gamma_j^{-1} (\bar \partial \gamma_j )\gamma_j^{-1}- \gamma_i \gamma_j^{-1} (\bar \partial \gamma_j 						)\gamma_j^{-1} \\ 
		&= 0, 
\end{split}  \end{equation*}
so we have associated to $\xi$ a holomorphic principal $G$-bundle. 
In Theorem \ref{thm:equiv} we have seen that two such bundles are holomorphically equivalent
if and only if $\xi_1=\xi_2*g$ for some $g \in \GS.$ Thus each orbit corresponds to a different 
equivalence class of holomorphic bundles. 
We have to show that every holomorphic principal $G$-bundle arises in this way. 
We have already seen that all bundles over $\Sigma$ are topologically trivial, 
so we let $P=M \times G.$ Let $\psi_j: (u,g) \mapsto (u, \gamma_j(u)g)$ denote the bundle charts. 
As the bundle is holomorphic, the transition functions are holomorphic, so 
$\bar \partial (\gamma_i \gamma_j^{-1}) = 0 = (\bar \partial \gamma_i )\gamma_j^{-1} - \gamma_i \gamma_j^{-1} (\bar \partial \gamma_j )\gamma_j^{-1}$ 
and therefore 
$\gamma_j^{-1}\bar \partial \gamma_j = \gamma_i^{-1} .\bar \partial \gamma_i$ 
on the intersection $U_i \cap U_j.$ Therefore the $\gamma_j^{-1} .\bar \partial \gamma_j$ 
define a $(0,1)$-form $\xi$ with values in $\g$ on $\Sigma,$ and obviously the bundle 
that we construct from $\xi$ as in the first part of our proof coincides with the given one.
\end{proof}

%
\subsection*{Acknowledgements}

The author wishes to thank Karl-Hermann Neeb for his continuing support, and Jimmie Lawson for encouraging me to submit this paper.

\LastPageEnding

\end{document}